\documentclass{amsart}
\usepackage{amssymb}
\usepackage[mathscr]{eucal}
\newtheorem{theorem}{Theorem}[section]

\newtheorem{lemma}[theorem]{Lemma}
\newtheorem{proposition}[theorem]{Proposition}
\theoremstyle{definition}
\newtheorem{definition}[theorem]{Definition}

\newtheorem{example}[theorem]{Example}

\theoremstyle{remark}

\begin{document}
\title{Lipschitzness of $^\ast$-homomorphisms between $C^\ast$-metric algebras}
\author{wei wu}
\thanks{The research was supported by the Shanghai Leading Academic Discipline Project (Project No.
B407), and the National Natural Science Foundation of China (Grant
No. 10671068).}
\address{Department of Mathematics, East China Normal University, Shanghai 200241, P.R. China}
\email{wwu@math.ecnu.edu.cn} 
\begin{abstract}A $C^\ast$-metric algebra consists of a unital
$C^\ast$-algebra and a Leibniz Lip-norm on the $C^\ast$-algebra.
We show that if the Lip-norms concerned are lower semicontinuous,
then any unital $^\ast$-homomorphism from a $C^\ast$-metric
algebra to another one is necessarily Lipschitz. It results that
the free product of two Lipschitz unital $^\ast$-homomorphisms
between $C^\ast$-metric algebras coming from $^\ast$-filtrations
is still a Lipschitz unital $^\ast$-homomorphism.
\end{abstract}
\maketitle

\section{Introduction}\label{se:1}

Originated in Kantorovi\v{c}'s work \cite{Kan,KanRub} and Connes'
observation \cite{Co1}, compact quantum metric space was
introduced by Rieffel as the noncommutative analogue of compact
metric space \cite{Ri1, Ri2}. It consists of an order unit space
and a Lip-norm on it, where the Lip-norm plays the role of the
usual Lipschitz seminorms for ordinary compact metric spaces
\cite{Ri5,OzRi,Ri6}. From this, a way to view ``matrix algebras
converge to the sphere" is established by defining the
Gromov-Hausdorff distance between compact quantum metric spaces
\cite{Ri3, Ri4}. See \cite{Ke1, Li1, Wu3, KeLi, Li2,Ri8} for
further discussion.

To precisely formulate the statements in the literature of
high-energy physics and string-theory, such as ``here are the
vector bundles over the matrix algebras that correspond to the
monopole bundles over the sphere" in \cite{Val,BaBaVaYd,GrRuSt}
etc., Rieffel found that algebraic conditions on the Lip-norm,
such as the Leibniz rule, play a crucial role \cite{Ri7, Ri8}.

At the investigation of the structure of unital $C^\ast$-algebras
with complexified Lip-norms, Kerr proposed a notion of dimension
for these $C^\ast$-algebras, along with two dynamical entropies.
In this environment, the Leibniz rule on the Lip-norm is also of
central importance \cite{Ke2}.

So unital $C^\ast$-algebras with Leibniz Lip-norms form an
important class of compact quantum metric spaces. For these class
of objects, it is essential that the maps preserve the given
algebraic structure and seminorms in same way. The morphisms
between two such objects were discussed in \cite{Ke2}: Lipschitz
unital $^\ast$-homomorphisms. It is natural to ask what conditions
are needed to guarantee that a unital $^\ast$-homomorphism from
one such object to another is Lipschitz. In this note we show

\begin{theorem}\label{th:11}
Let $(A_1,L_1)$ and $(A_2,L_2)$ be two $C^\ast$-metric algebras
with lower semicontinuous Lip-norms. If $T$ is a unital
$^\ast$-homomorphism from $(A_1,L_1)$ into $(A_2,L_2)$, then $T$
is Lipschitz.
\end{theorem}
\noindent (Precise definitions are given in Section \ref{se:2}.)

In fact, any Lip-norm on a unital $C^\ast$-algebra may be replaced
by a lower semicontinuous Lip-norm according to theorem 4.1 in
\cite{Ri2}. Therefore, unital $C^\ast$-algebras with lower
semicontinuous Lip-norms are a class of $C^\ast$-algebras that
deserve to be studied further.

We recall in Section \ref{se:2} some terminology and some typical
constructions of unital $C^\ast$-algebras with lower
semicontinuous Leibniz Lip-norms. We prove Theorem \ref{th:11} in
Section \ref{se:3}, and as an application we show, in Section
\ref{se:4}, that the free product of two Lipschitz unital
$^\ast$-homomorphisms between $C^\ast$-metric algebras coming from
$^\ast$-filtrations is still a Lipschitz unital
$^\ast$-homomorphism.

Acknowledgements. I am grateful to Hanfeng Li for comments.

\section{$C^\ast$-metric algebras}\label{se:2}

We recall \cite{Ri8} that a {\it Lip-norm} on a unital
$C^\ast$-algebra $A$ is a seminorm $L$ on $A$, possibly taking the
value $+\infty$, such that
\begin{enumerate}
\item $L(1_A)=0$, where $1_A$ is the identity of $A$,
\item $L(a^\ast)=L(a)$ for all $a\in A$,
\item\label{en:1} the metric $\rho_L$ on $S(A)$ given by
\[\rho_L(\varphi,\psi)=\sup\{|\varphi(a)-\psi(a)|:L(a)\leq 1\}\]
induces the weak$^\ast$-topology on $S(A)$.
\end{enumerate}

When $L$ is a Lip-norm on a unital $C^\ast$-algebra $A$, the
metric space $(S(A),\rho_L)$ is compact. So the diameter
$\mathrm{diam}(A,L)$ of $(A,L)$, given by the diameter of $S(A)$
with respect to $\rho_L$, is finite. If $L(a)=0$ and the condition
\ref{en:1} above holds, then $L(na)=0$ for any $n\in\mathbb{N}$,
and so $|\varphi(a)-\psi(a)|\leq\frac1n\rho_L(\varphi,\psi)$ for
all $\varphi,\psi\in S(A)$ and $n\in\mathbb{N}$. Hence
$\varphi(a)=\psi(a)$ for all $\varphi,\psi\in S(A)$. Therefore,
$a\in\mathbb{C}1_A$. Let $\varphi,\psi\in S(A)$ with
$\varphi\neq\psi$. Since the weak$^\ast$-topology on $S(A)$ is
Hausdorff, $\rho_L(\varphi,\psi)>0$ if the metric $\rho_L$ gives
$S(A)$ the weak$^\ast$-topology. So there is an $a\in A$ with
$L(a)\leq 1$ such that $\varphi(a)\neq\psi(a)$. This indicates
that $L$ separates the state space $S(A)$ of $A$.

A seminorm $L$ on a unital $C^\ast$-algebra $A$ is said to {\it
lower semicontinuous} if for one $r\in\mathbb{R}_{>0}$, hence for
all $r>0$, the set $\{a\in A: L(a)\leq r\}$ is norm-closed in $A$.
Given a seminorm on a unital $C^\ast$-algebra $A$, we denote the
set $\{a\in A:L(a)<\infty\}$ of Lipschitz elements in $A$ by
$\mathcal{A}$. The seminorm $L$ is said to be {\it Leibniz} if it
satisfies the Leibniz rule
\[L(ab)\leq L(a)\|b\|+\|a\|L(b)\]
for all $a,b\in\mathcal{A}$.

Now we define

\begin{definition}\label{de:21}
A $C^\ast$-{\it metric algebra} is a pair $(A,L)$ consisting of a
unital $C^\ast$-algebra $A$ and a Leibniz Lip-norm $L$ on $A$.
\end{definition}

Basically, compact metric spaces provide the prototype of
$C^\ast$-metric algebras.

\begin{example}Let $(X,\rho)$ be an ordinary compact metric space,
and let $L_{\rho}$ denote the Lipschitz seminorm on $C(X)$. Then
$(C(X),L_{\rho})$ is a $C^\ast$-metric algebra and $L_{\rho}$ is
lower semicontinuous.
\end{example}

One important class of examples of $C^\ast$-metric algebras comes
from ergodic actions of a compact group on a unital
$C^\ast$-algebra \cite{Ri1,Ri8}.

\begin{example}Let $A$ be a unital $C^\ast$-algebra, let $G$ be a compact group,
and let $\alpha$ be a strongly continuous ergodic action of $G$ on
$A$ by automorphisms. We assume that $G$ is equipped with a length
function $\ell$. Then the seminorm
\[L(a)=\sup\left\{\|\alpha_x(a)-a\|/\ell(x):x\in G\backslash\{e\}\right\},\ \ \ a\in A,\]
is a lower semicontinuous Leibniz Lip-norm on $A$.
\end{example}

In \cite{Li3}, a generalization to ergodic actions of co-amenable
compact quantum group was given. By viewing an ergodic action of a
compact group on a unital $C^\ast$-algebra as the translation
action of the group on a noncommutative homogeneous space of it,
Hanfeng Li recently extended the construction above to the locally
compact groups which satisfy certain conditions \cite{Li4}.

\begin{example}
Given a spectral triple $(A,H,D)$, where $(H,D)$ is an unbounded
Fredholm module over a unital $C^\ast$-algebra $A$ \cite{Co2,Co3},
we can define a seminorm, $L$, on $A$ by
\[L(a)=\|[D,a]\|.\]
Let $\sigma$ be a state of $A$. Then $L$ is a Lip-norm if and only
the set $\{a\in A:L(a)\leq 1,\sigma(a)=0\}$ is a
norm-totally-bounded subset of $A$ \cite{OzRi}. In this case,
$(A,L)$ is a $C^\ast$-metric algebra with lower semicontinuous
Lip-norm.
\end{example}

In fact, every $C^\ast$-metric algebra can be obtained from the
Dirac operator approach \cite{Ri6}. Finally, we make:

\begin{definition}\label{de:22}
Let $(A_1,L_1)$ and $(A_2,L_2)$ be two $C^\ast$-metric algebras.
By a $^\ast$-{\it homomorphism from} $(A_1,L_1)$ {\it to}
$(A_2,L_2)$ we mean a $\ast$-homomorphism $T$ from
$C^\ast$-algebra $A_1$ to $C^\ast$-algebra $A_2$ such that
$T(\mathcal{A}_1)\subseteq\mathcal{A}_2$. A map $\alpha:
A_1\rightarrow A_2$ is said to be {\it Lipschitz} if there is a
constant $\mu\geq 0$ such that
\[L_2(\alpha(a))\leq\mu L_1(a)\]
for all $a\in\mathcal{A}_1$.
\end{definition}

\section{$^\ast$-homomorphisms between $C^\ast$-metric algebras}\label{se:3}

Let $A$ be a unital $C^\ast$-algebra with a {\it Lipschitz
seminorm} $L$, a seminorm which may take value $+\infty$, has its
null-space the scalar multiples of the identity, and satisfies
$L(a^\ast)=L(a)$ for all $a\in A$. We can define a metric $\rho_L$
on the state space $S(A)$ of $A$ by
\[\rho_L(\varphi,\psi)=\mathrm{sup}\{|\varphi(a)-\psi(a)|:L(a)\leq 1\},\ \ \ \varphi,\psi\in S(A).\]
Usually this metric may take value $+\infty$. From this metric we
can define a seminorm $L_{\rho_L}$ on $A$ by the formula
\[L_{\rho_L}(a)=\mathrm{sup}\left\{\frac{|\varphi(a)-\psi(a)|}{\rho_L(\varphi,\psi)}:\varphi\neq\psi,\varphi,\psi\in S(A)\right\},\ \ \ a\in A.\]
This is still a Lipschitz seminorm.

The dual Banach space to $A/(\mathbb{C}1_A)$ for the quotient norm
$\|\cdot\|^\sim$ is just the subspace $A^{\prime 0}$ of the dual
space $A^\prime$ consisting of those $f\in A^\prime$ such that
$f(1_A)=0$. Similar to \cite{Ri2}, we define
\[L^\prime(f)=\mathrm{sup}\{|f(a)|:L(a)\leq 1,a\in A\},\ \ \ f\in A^{\prime}.\]
By proposition 3.7 of \cite{Pa} and proposition 4.2 of \cite{Wu1},
for any $g\in A^{\prime 0}$ with $\|g\|\leq 2$ there are
$\varphi_1,\varphi_2,\varphi_3,\varphi_4\in S(A)$ such that
$g=(\varphi_1-\varphi_2)+i(\varphi_3-\varphi_4)$. So
\begin{align*}
L^\prime(g)&=\mathrm{sup}\{|g(a)|:L(a)\leq 1,a\in A\}\\
&=\mathrm{sup}\{|(\varphi_1(a)-\varphi_2(a))+i(\varphi_3(a)-\varphi_4(a))|:L(a)\leq
1,a\in A\}\\
&\leq\mathrm{sup}\{|\varphi_1(a)-\varphi_2(a)|:L(a)\leq 1,a\in
A\}\\
&+\mathrm{sup}\{|\varphi_3(a)-\varphi_4(a)|:L(a)\leq 1,a\in
A\}\\
&=\rho_L(\varphi_1,\varphi_2)+\rho_L(\varphi_3,\varphi_4),
\end{align*}
and
\begin{align*}
\rho_L(\varphi_1,\varphi_2)+\rho_L(\varphi_3,\varphi_4)
&=\mathrm{sup}\{|\varphi_1(a)-\varphi_2(a)|:L(a)\leq 1,a\in
A\}\\
&+\mathrm{sup}\{|\varphi_3(a)-\varphi_4(a)|:L(a)\leq 1,a\in
A\}\\
&=\mathrm{sup}\{|g(a)+g^\ast(a)|/2:L(a)\leq 1,a\in
A\}\\
&+\mathrm{sup}\{|g(a)-g^\ast(a)|/2:L(a)\leq 1,a\in
A\}\\
&\leq 2\mathrm{sup}\{|g(a)|:L(a)\leq 1,a\in A\}\\
&=2L^\prime(g).
\end{align*}

Denote
\[\mathcal{L}_1=\{a\in A:L(a)\leq 1\}.\]
Then $\mathcal{L}_1$ is convex and balanced, and the bipolar
theorem says that
$\overline{\mathcal{L}}_1^{\|\cdot\|}=\mathcal{L}_1^{\circ\circ}$.
Now suppose that $L$ is lower semicontinuous. Then we get
\[\mathcal{L}_1=\mathcal{L}_1^{\circ\circ}.\]
If $L_{\rho_L}(a)\leq 1$, then
$|\varphi(a)-\psi(a)|\leq\rho_L(\varphi,\psi)$ for all
$\varphi,\psi\in S(A)$. So $|g(a)|\leq 2L^\prime(g)$ for all $g\in
A^{\prime0}$ with $\|g\|\leq 2$, and this also implies that
$|g(a)|\leq 2L^\prime(g)$ for all $g\in A^{\prime0}$. Thus
$|g(a/2)|\leq 1$ for all $g\in A^{\prime}$ with $L^\prime(g)\leq
1$. By definition, we obtain that
$a/2\in\mathcal{L}_1^{\circ\circ}=\mathcal{L}_1$. So
\[L(a)\leq 2L_{\rho_L}(a).\]
If $L(a)\leq 1$, then $a\in\mathcal{L}_1^{\circ\circ}$. So
$|f(a)|\leq 1$ for all $f\in A^{\prime}$ with $L^\prime(f)\leq 1$,
and hence $|f(a)|\leq L^\prime(f)$ for all $f\in A^{\prime0}$. In
particular, we have that $|\varphi(a)-\psi(a)|\leq
L^\prime(\varphi-\psi)=\rho_L(\varphi,\psi)$ for $\varphi,\psi\in
S(A)$. Thus $L_{\rho_L}(a)\leq 1$. So
\[L_{\rho_L}(a)\leq L(a).\]
Therefore, we have

\begin{lemma}\label{le:30}
Let $A$ be a unital $C^\ast$-algebra with a lower semicontinuous
Lipschitz seminorm $L$. Then
\[L_{\rho_L}(a)\leq L(a)\leq 2L_{\rho_L}(a),\]
for all $a\in A$.
\end{lemma}

Suppose $A_1$ and $A_2$ are two unital $C^\ast$-algebras with
$C^\ast$-norms $\|\cdot\|_1$ and $\|\cdot\|_2$, respectively. Let
$L_1$ and $L_2$ be Lipschitz seminorms on $A_1$ and $A_2$,
respectively, which are Leibniz. On $\mathcal{A}_1$ and
$\mathcal{A}_2$, we can define new norms by
\[M_1(a)=\|a\|_1+L_1(a),\ \ \ a\in \mathcal{A}_1,\]
\[M_2(b)=\|b\|_2+L_2(b),\ \ \ b\in \mathcal{A}_2,\]
respectively. Then $(\mathcal{A}_1, M_1)$ and $(\mathcal{A}_2,
M_2)$ are normed algebras since $L_1$ and $L_2$ are Leibniz. We
denote by $\mathfrak{A}_1$ and $\mathfrak{A}_2$ the completions of
$\mathcal{A}_1$ and $\mathcal{A}_2$ with respect to the norms
$M_1$ and $M_2$, respectively.

Let $\{a_n\}$ be a sequence in $\mathcal{A}_1$ such that
$\sum_{n=1}^\infty M_1(a_n)$ converges. Then both
$\sum_{n=1}^\infty \|a_n\|_1$ and $\sum_{n=1}^\infty L_1(a_n)$ are
convergent. Since $A_1$ is complete in the norm $\|\cdot\|_1$,
there is an $a\in A_1$ such that $a=\sum_{n=1}^\infty a_n$ in the
norm $\|\cdot\|_1$. For any $k\in\mathbb{N}$, we have
\[L_1\left(\sum_{n=1}^k a_n\right)\leq\sum_{n=1}^k L_1(a_n)\leq\sum_{n=1}^\infty L_1(a_n).\]
Now assume that $L_1$ is lower semicontinuous. Then we have
\[L_1(a)\leq\sum_{n=1}^\infty L_1(a_n)<\infty.\]
Thus $a\in\mathcal{A}_1$.
Similarly, we have that $L_1\left(\sum_{n=k+1}^\infty
a_n\right)\leq\sum_{n=k+1}^\infty L_1(a_n)$ for $k\in\mathbb{N}$.
So
\begin{align*}
M_1\left(a-\sum_{n=1}^k a_n\right)&=\left\|a-\sum_{n=1}^k
a_n\right\|_1+L_1\left(a-\sum_{n=1}^k
a_n\right)\\
&\leq\sum_{n=k+1}^\infty \|a_n\|_1+\sum_{n=k+1}^\infty
L_1(a_n).\end{align*} Thus
$\lim_{k\rightarrow\infty}M_1\left(a-\sum_{n=1}^k a_n\right)=0$,
i.e., $\sum_{n=1}^\infty a_n$ converges to $a$ in the norm $M_1$.
So $\mathcal{A}_1$ is complete in the norm $M_1$. Therefore,
$\frak{A}_1=\mathcal{A}_1$. Similarly, we have that
$\frak{A}_2=\mathcal{A}_2$ if $L_2$ is lower semicontinuous.
Moreover, in this case we have

\begin{proposition}\label{pro:31}
If $\alpha$ is a $^\ast$-homomorphism from $C^\ast$-algebra $A_1$
into $C^\ast$-algebra $A_2$ such that
$\alpha(\mathcal{A}_1)\subseteq\mathcal{A}_2$, then there is a
constant $\lambda>0$ such that
\[M_2(\alpha(a))\leq\lambda M_1(a)\]
for all $a\in\mathcal{A}_1$.
\end{proposition}

\begin{proof}
Since both the norms $\|\cdot\|_1$ and $\|\cdot\|_2$ and the
seminorms $L_1$ and $L_2$ have an isometric involution, the
involutions on $\mathcal{A}_1$ and $\mathcal{A}_2$ are isometric
with the norms $M_1$ and $M_2$, respectively. And so
$(\mathcal{A}_1, M_1)$ and $(\mathcal{A}_2, M_2)$ are Banach
$^\ast$-algebras with isometric involutions. That $\alpha$ is a
$^\ast$-homomorphism implies that $\alpha(\mathcal{A}_1)$ is a
$^\ast$-subalgebra of $\mathcal{A}_2$. Also from the relation
$\mathcal{A}_2\subseteq {A}_2$, we get the restriction of the
$C^\ast$-norms on $A_2$ to $\mathcal{A}_2$, and it makes
$\mathcal{A}_2$ into a pre-$C^\ast$-algebra \cite{Wu4}. This shows
that $(\mathcal{A}_2, M_2)$ is an $A^\ast$-algebra. By theorem
23.11 in \cite{DoBe}, $\alpha$ is continuous as a map from
$(\mathcal{A}_1, M_1)$ into $(\mathcal{A}_2, M_2)$. Thus, there is
a constant $\lambda>0$ such that $M_2(\alpha(a))\leq\lambda
M_1(a)$ for all $a\in\mathcal{A}_1$.
\end{proof}

Assume that $T$ is a unital positive linear map from
$C^\ast$-algebra $A_1$ into $C^\ast$-algebra $A_2$. Let $S(A_1)$
and $S(A_2)$ denote the state spaces of $A_1$ and $A_2$,
respectively. For any $\varphi\in S(A_2)$, we define
\[\sigma_T(\varphi)(a)=\varphi(Ta),\ \ \ a\in A_1.\]
Then $\sigma_T$ is an affine map from $S(A_2)$ into $S(A_1)$.
Moreover, if $T$ maps $A_1$ onto $A_2$, then $\sigma_T$ is
injective; if $\sigma_T$ is a surjection, then $T$ is an
injection.

\begin{lemma}\label{le:31}
Let $(A_1,L_1)$ and $(A_2,L_2)$ be two $C^\ast$-metric algebras
with lower semicontonuous Lip-norms. If $T$ is a unital
$^\ast$-homomorphism from $(A_1,L_1)$ into $(A_2,L_2)$, then there
is a Lipschitz affine map $\sigma$ from metric space
$(S(A_2),\rho_{L_2})$ into metric space $(S(A_1),\rho_{L_1})$ such
that
\[\varphi(Ta)=\sigma(\varphi)(a), \ \ \ a\in A_1, \varphi\in S(A_2).\]
\end{lemma}

\begin{proof}
Suppose $T$ is a unital $^\ast$-homomorphism from $A_1$ into
$A_2$. Set $\sigma=\sigma_T$. For any $\varphi_0\in S(A_1)$ and
$b\in A_1$ with $L_1(b)\leq 1$, define
\[f_{\varphi_0,b}(\varphi)=\varphi(b)-\varphi_0(b),\ \ \ \varphi\in S(A_1).\]
Clearly $f_{\varphi_0,b}$ is a continuous affine function on
$S(A_1)$. To each $a\in A_1$, we get a function $\hat{a}$ on
$S(A_1)$ defined by $\hat{a}(\varphi)=\varphi(a)$. By the basic
representation theorem of Kadison \cite{Ka} (see theorem II.1.8 of
\cite{Al}), there is a self-adjoint $a_{\varphi_0,b}\in A_1$ such
that $f_{\varphi_0,b}=\hat{a}_{\varphi_0,b}$ for every
$b=b^\ast\in A_1$.

For any $\varphi\in S(A_1)$ and $b=b^\ast\in A_1$ with $L_1(b)\leq
1$, we have
\[|\hat{a}_{\varphi_0,b}(\varphi)|=|\varphi(b)-\varphi_0(b)|\leq\rho_{L_1}(\varphi,\varphi_0)\leq\mathrm{diam}(A_1,L_1).\]
So $\|a_{\varphi_0,b}\|_1\leq\mathrm{diam}(A_1,L_1)$ by
proposition II.1.7 of \cite{Al}. And since $L_1$ is lower
semicontinuous, we get
\begin{align*}
L_1(a_{\varphi_0,b})&\leq 2L_{\rho_{L_1}}(a_{\varphi_0,b})\\
&=2\mathrm{sup}\left\{\frac{|\varphi_1(a_{\varphi_0,b})-\varphi_2(a_{\varphi_0,b})|}{\rho_{L_1(\varphi_1,\varphi_2)}}:\varphi_1\neq\varphi_2,\varphi_1,\varphi_2\in
S(A_1)\right\}\\
&=2\mathrm{sup}\left\{\frac{|\hat{a}_{\varphi_0,b}(\varphi_1)-\hat{a}_{\varphi_0,b}(\varphi_2)|}{\rho_{L_1(\varphi_1,\varphi_2)}}:\varphi_1\neq\varphi_2,\varphi_1,\varphi_2\in
S(A_1)\right\}\\
&=2\mathrm{sup}\left\{\frac{|\varphi_1(b)-\varphi_2(b)|}{\rho_{L_1(\varphi_1,\varphi_2)}}:\varphi_1\neq\varphi_2,\varphi_1,\varphi_2\in
S(A_1)\right\}\\
&=2L_{\rho_{L_1}}(b)\leq 2L_1(b)\leq 2
\end{align*}
by Lemma \ref{le:30}. Thus
\[M_1(a_{\varphi_0,b})\leq\mathrm{diam}(A_1,L_1)+2.\]
So the set $\{a_{\varphi_0,b}: b=b^\ast\in A_1, L_1(b)\leq 1\}$ is
bounded in $A_1$ with respect to the norm $M_1$. By Proposition
\ref{pro:31}, the set $\{Ta_{\varphi_0,b}: b=b^\ast\in A_1,
L_1(b)\leq 1\}$ is bounded in $A_2$ under the norm $M_2$. Then
there is a constant $K>0$ such that $M_2(Ta_{\varphi_0,b})\leq K$
for all $b=b^\ast\in A_1$ with $L_1(b)\leq 1$.

Now for any $\psi_1,\psi_2\in S(A_2)$ and any $b=b_1+ib_2\in A_1$
with $L_1(b)\leq 1$ and $b_i=b_i^\ast\in A_1$ for $i\in\{1,2\}$,
we have that $L_1(b_1)\leq 1$ and $L_1(b_2)\leq 1$, and so
\begin{align*}
|\sigma(\psi_1)(b_i)-\sigma(\psi_2)(b_i)|
&=|[\sigma(\psi_1)(b_i)-\varphi_0(b_i)]-[\sigma(\psi_2)(b_i)-\varphi_0(b_i)]|\\
&=|\hat{a}_{\varphi_0,b_i}(\sigma(\psi_1))-\hat{a}_{\varphi_0,b_i}(\sigma(\psi_2))|\\
&=|\hat{Ta}_{\varphi_0,b_i}(\psi_1)-\hat{Ta}_{\varphi_0,b_i}(\psi_2)|\\
&\leq
L_{\rho_{L_2}}(Ta_{\varphi_0,b_i})\rho_{L_2}(\psi_1,\psi_2)\\
&\leq L_2(Ta_{\varphi_0,b_i})\rho_{L_2}(\psi_1,\psi_2)\\
&\leq M_2(Ta_{\varphi_0,b_i})\rho_{L_2}(\psi_1,\psi_2)\\
&\leq K\rho_{L_2}(\psi_1,\psi_2),
\end{align*}
by Lemma \ref{le:30}. From this we obtain
\begin{align*}
|\sigma(\psi_1)(b)-\sigma(\psi_2)(b)|
&\leq|\sigma(\psi_1)(b_1)-\sigma(\psi_2)(b_1)|+|\sigma(\psi_1)(b_2)-\sigma(\psi_2)(b_2)|\\
&\leq 2K\rho_{L_2}(\psi_1,\psi_2).\end{align*} Therefore,
$\rho_{L_1}(\sigma(\psi_1),\sigma(\psi_2))\leq
2K\rho_{L_2}(\psi_1,\psi_2)$. Hence $\sigma$ is Lipschitz.
\end{proof}

We are ready to prove Theorem \ref{th:11}.

\medskip
\noindent
{\it Proof of Theorem \ref{th:11}} By Lemma \ref{le:31},
there exists an affine map $\sigma$ from metric space
$(S(A_2),\rho_{L_2})$ into metric space $(S(A_1),\rho_{L_1})$ such
that
\[\varphi(Ta)=\sigma(\varphi)(a), \ \ \ a\in A_1, \varphi\in S(A_2),\]
and
\[\rho_{L_1}(\sigma(\psi_1),\sigma(\psi_2))\leq K\rho_{L_2}(\psi_1,\psi_2), \ \ \ \psi_1,\psi_2\in S(A_2),\]
for some positive constant $K$. So for any $a\in \mathcal{A}_1$,
we have
\begin{align*}
L_2(Ta)&\leq 2L_{\rho_{L_2}}(Ta)\\
&=2\mathrm{sup}\left\{\frac{|\varphi(Ta)-\psi(Ta)|}{\rho_{L_2}(\varphi,\psi)}:\varphi\neq\psi,\varphi,\psi\in
S(A_2)\right\}\\
&=2\mathrm{sup}\left\{\frac{|\sigma(\varphi)(a)-\sigma(\psi)(a)|}{\rho_{L_2}(\varphi,\psi)}:\varphi\neq\psi,\varphi,\psi\in
S(A_2)\right\}\\
&=2\mathrm{sup}\left\{\frac{|\sigma(\varphi)(a)-\sigma(\psi)(a)|}{\rho_{L_1}(\sigma(\varphi),\sigma(\psi))}\cdot
\frac{\rho_{L_1}(\sigma(\varphi),\sigma(\psi))}{\rho_{L_2}(\varphi,\psi)}:\varphi\neq\psi,\varphi,\psi\in
S(A_2)\right\}\\
&\leq 2KL_1(a)
\end{align*}
by Lemma \ref{le:30}. Therefore, $T$ is Lipschitz.$\hfill\square$

\section{an application}\label{se:4}

Let $A$ be a unital $C^\ast$-algebra. A $^\ast$-{\it filtration}
$\{A_n\}$ of $A$ is a sequence of finite-dimensional subspaces
which satisfy
\begin{enumerate}
\item $A_0=\mathbb{C}1_A$, $A^\ast_n=A_n$,
\item $A_m\subset A_n$ if $m<n$,
\item $A_mA_n\subseteq A_{m+n}$,
\item $A=\overline{\cup^\infty_{n=0}A_n}$
\end{enumerate}
\cite{Vo3}. Given a faithful state $\sigma$ on $A$. Let
$(\pi,H,\xi)$ be the faithful GNS representation of $(A,\sigma)$.
We identify $A$ with the corresponding linear space of vectors in
$H$. We let $\|\cdot\|$ denote the operator norm of $\pi(A)$, and
$\|\cdot\|_2$ denote the vector norm of $A$. Viewing each $A_n$ as
a finite-dimensional subspace of $H$, we obtain a filtration of
$H$ in the sense that $\{A_n\}$ is an increasing sequence of
finite-dimensional subspaces of $H$ and $\cup^\infty_{n=0}A_n$ is
dense in $H$ \cite{Vo3}. Let $Q_n$ denote the orthogonal
projection of $H$ onto $A_n$. Set $P_n=Q_n-Q_{n-1}$ for $n\geq 1$,
and $P_0=Q_0$. We define
\[D=\sum^\infty_{n=1}nP_n.\]
Then $D$ is a unbounded linear operator on $H$ with domain
$\mathcal{A}=\cup^\infty_{n=0}A_n$. The linear functional
$f(h)=\langle Dh,k\rangle$ on $\mathcal{A}$ is bounded if and only
if $k\in\mathcal{A}$. It is also clear that $\langle
Dh,k\rangle=\langle h,Dk\rangle$ for all $h,k\in\mathcal{A}$. So
$D$ is self-adjoint \cite{Con}. By lemma 1.1 of \cite{OzRi},
$\mathcal{A}=\{a\in A:[D,\pi(a)]\mbox{ is bounded}\}$. Furthermore
$D$ has a finite-dimensional kernel, and the inverse $D^{-1}$
(defined on the orthogonal complement of $D$'s kernel) is compact.
So $(\mathcal{A},H,D)$ is a spectral triple
\cite{Co1,Co2,Co3,Lan,Va}. Using this spectral triple, we can
define a seminorm $L$ on $A$, possible taking the value $+\infty$,
by
\[L(a)=\|[D,\pi(a)]\|,\ \ \ \ a\in A.\]
For all $a,b\in A$, we have
\begin{align*}
L(ab)&=\|[D,\pi(ab)]\|=\|[D,\pi(a)]\pi(b)+\pi(a)[D,\pi(b)]\|\\
&\leq\|[D,\pi(a)]\|\|\pi(b)\|+\|\pi(a)\|\|[D,\pi(b)]\|\\
&=L(a)\|b\|+\|a\|L(b). \end{align*} So $L$ is Leibniz. By main
theorem 1.2 in \cite{OzRi}, we have

\begin{proposition}\label{pr:41}
If there is a constant $C$ such that
\[\|P_m\pi(P_k(a))P_n\|\leq C\|P_k(a)\|_2\]
for all $a\in\mathcal{A}$ and $m,n,k\in\mathbb{N}$, then $(A,L)$
is a $C^\ast$-metric algebra with lower semicontinuous Lip-norm.
\end{proposition}

A condition of this kind is called a {\it Haagerup-type condition}
with constant $C$. In this case, we will call $(A,L)$ the
$C^\ast$-{\it metric algebra coming from} $(\{A_n\},\sigma)$ (with
constant $C$).

Suppose that $(A^1,L_1)$ and $(A^2,L_2)$ are two $C^\ast$-metric
algebras coming from $(\{A^1_n\},\sigma_1)$ and
$(\{A^2_n\},\sigma_2)$, respectively, and both with constant $C$.
Let
\[(A,\sigma)=(A^1,\sigma_1)\ast(A^2,\sigma_2)\]
be the reduced free-product $C^\ast$-algebra with the faithful
state $\sigma$\cite{Vo1,VoDyNi,BlDy}. We define a
$^\ast$-filtration $\{A_n\}$ on $(A,\sigma)$ by setting $A_n$ to
be the linear span of all products
$A^{\imath_1}_{n_1}\cdot\cdots\cdot A^{\imath_k}_{n_k}$ with each
$\imath_j\in\{1,2\}$, with $\imath_j\neq \imath_{j+1}$ for $1\leq
j\leq k-1$, and with $\sum^{k}_{j=1}n_j\leq n$. We let
$(\pi_\imath,H^\imath,\xi_\imath)$ denote the faithful GNS
representation of $(A^\imath,\sigma_{\imath})$ for
$\imath\in\{1,2\}$, and we let $(\pi,H,\xi)$ denote the faithful
GNS representation of $(A,\sigma)$. We let $\{P_n\}$ be the family
of mutually orthogonal projections corresponding to the filtration
$\{A_n\}$ as above. Let $D$ be the Dirac operators coming from
$(A,\{A_n\},\sigma)$. Then
\[D=\sum_{n=1}^{\infty}nP_n,\]
and the corresponding lower semicontinuous seminorm on $A$ is
\[L(a)=\|[D,\pi(a)]\|,\ \ a\in A.\]
By theorem 6.1 of \cite{OzRi}, we obtain

\begin{proposition}\label{pr:42}
If $\sigma_1$ and $\sigma_2$ are traces, then
$(A,L)$is a $C^\ast$-metric algebra coming from $(\{A_n\},\sigma)$
with constant $\sqrt{5}C$, the reduced $C^\ast$-norm on $A$ is the
operator norm for the GNS representation for $\sigma$ on $H$, and
$(H,\xi)=(H^1,\xi_1)\ast(H^2,\xi_2)$.
\end{proposition}

From the discussion above, we see that
$\mathcal{A}=\cup_{n=0}^\infty
A_n=\mathcal{A}_1\ast\mathcal{A}_2$, where
$\mathcal{A}_1\ast\mathcal{A}_2$ is the algebraic free-product of
$\mathcal{A}_1$ and $\mathcal{A}_2$ with its evident involution.

\begin{theorem}\label{th:41}
Let $(A_\imath,L_\imath)$ and
$(\tilde{A}_\imath,\tilde{L}_\imath)$ be $C^\ast$-metric algebras
coming from $(\{A_\imath\},\sigma_\imath)$ and
$(\{\tilde{A}_\imath\},\tilde{\sigma}_\imath)$, respectively, with
the faithful tracial states $\sigma_{\imath}$ and
$\tilde{\sigma}_{\imath}$ and constant $C$ for all $\imath=1,2$.
Denote
\[(A,\sigma)=(A_1,\sigma_1)\ast(A_2,\sigma_2),\]
\[(\tilde{A},\tilde{\sigma})=(\tilde{A}_1,\tilde{\sigma}_1)\ast(\tilde{A}_2,\tilde{\sigma}_2).\]
Suppose that $\alpha_1: A_1\rightarrow\tilde{A}_1$ and $\alpha_2:
A_2\rightarrow\tilde{A}_2$ are two Lipschitz unital
$\ast$-homomorphisms such that
$\tilde{\sigma}_\imath\circ\alpha_\imath=\sigma_\imath$. Then
there is a unique Lipschitz unital $\ast$-homomorphism
$\alpha:A\rightarrow\tilde{A}$ such that for every
$\imath\in\{1,2\}$ the diagram
\[\begin{array}{ccccc}
  \hspace{0.5cm}A_\imath & \stackrel{\lambda_\imath}{\longrightarrow} &\hspace{0.3cm}A& & \\
  \alpha_\imath\downarrow &  & \alpha\downarrow& \hspace{0.3cm}\searrow\sigma \\
  \hspace{0.5cm}\tilde{A}_\imath & \stackrel{\tilde{\lambda}_\imath}{\longrightarrow} & \hspace{0.3cm}\tilde{A}&\stackrel{\tilde{\sigma}}{\longrightarrow} & \mathbb{\mathbb{C}}
\end{array}\]
commutes, where $\lambda_\imath$ and $\tilde{\lambda}_\imath$ are
the maps arising from the free product construction.
\end{theorem}

\begin{proof}
Consider the unital $^\ast$-homomorphisms
$\beta_1=\tilde{\lambda}_1\circ\alpha_1:A_1\rightarrow\tilde{A}$
and
$\beta_2=\tilde{\lambda}_2\circ\alpha_2:A_2\rightarrow\tilde{A}$.
We have the following diagrams commute:
\[\begin{array}{ccccc}
     &                      &      \mathbb{C}\hspace{0.2cm}     &                     &          \\
     &\nearrow_{\sigma_1} & \uparrow_{\tilde{\sigma}}       &\hspace{0.3cm}\nwarrow_{\sigma_2}&     \\
A_1&\stackrel{\beta_1}{\longrightarrow}&\tilde{A}\hspace{0.2cm}&\stackrel{\beta_2}{\longleftarrow}&A_2
\end{array}\]
for $\imath\in\{1,2\}$. By lemma 1.3 of \cite{DyRo}, there is a
$^\ast$-homomorphism $\alpha:A\rightarrow\tilde{A}$ such that
$\alpha\circ\lambda_{\imath}=\beta_{\imath}$ and
$\tilde{\sigma}\circ\alpha=\sigma$. In particular, we have that
$\alpha\circ\lambda_{\imath}=\tilde{\lambda}_\imath\circ\alpha_\imath$
and
$\alpha(1_{A})=\alpha(\lambda_1(1_{A_1}))=\beta_1(1_{A_1})=(\tilde{\lambda}_1\circ\alpha_1)(1_{A_1})=1_{\tilde{A}}$.

Since $A$ is generated by $\lambda_1(A_1)\cup\lambda_2(A_2)$, it
is clear that $\alpha$ will be unique if it exists. It is also
clear that $\alpha(\mathcal{A})\subseteq\tilde{\mathcal{A}}$, and
so the Lipschitzness of $\alpha$ follows from Proposition
\ref{pr:41}, Proposition \ref{pr:42} and Theorem \ref{th:11}.
\end{proof}

\bibliographystyle{amsplain}

\begin{thebibliography}{20}

\bibitem{Al} E. M. Alfsen, \textit{Compact convex sets and boundary
integrals}, Springer-Verlag, Berlin, New York, 1971.

\bibitem{BaBaVaYd} S. Baez, A. P. Balachandran, S. Vaidya, B. Ydri, \textit{Monopoles and solitons in fuzzy physics}, Comm. Math. Phys., \textbf{208}(2000),
no. 3, 787-798. arXiv:hep-th/9811169.

\bibitem{BlDy} E. F. Blanchard, K. J. Dykema, \textit{Embeddings of reduced free products of operator algebras},
Pacific J. Math., \textbf{199}(2001), no. 1, 1-19.
arXiv:math.OA/9911012.

\bibitem{Co1} A. Connes, \textit{Compact metric spaces, Fredholm modules, and
hyperfiniteness}, Ergodic Theory Dynam. Systems, \textbf{9}(1989),
no. 2, 207-220.

\bibitem{Co2} A. Connes, \textit{Noncommutative geometry}, Academic
Press, Inc., San Diego, CA, 1994.

\bibitem{Co3} A. Connes, \textit{Noncommutative geometry and
reality}, J. Math. Phys., \textbf{36}(1995), no. 11, 6194-6231.

\bibitem{Con} J. B. Conway, \textit{A course in functional
analysis}, Springer-Verlag, New York, second edition, 1990.

\bibitem{DoBe} R. S. Doran, V. A. Belfi, \textit{Characterizations of $C^\ast$-algebras.
The Gelfand-Naimark theorems}, Monographs and Textbooks in Pure
and Applied Mathematics, 101. Marcel Dekker, Inc., New York, 1986.

\bibitem{DyRo} K. J. Dykema, M. R\o rdam, \textit{Projections in free product $C^\ast$-algebras}, Geom. Funct. Anal., \textbf{8}(1998),
1-16. arXiv:funct-an/9702016.

\bibitem{GrRuSt} Harald Grosse, Christian W. Rupp, Alexander Strohmaier, \textit{Fuzzy line bundles, the Chern character and topological charges over
the fuzzy sphere}, J. Geom. Phys., \textbf{42}(2002), no. 1-2,
54-63. arXiv:math-ph/0105033.

\bibitem{Ka} R. V. Kadison, \textit{A representation theory for commutative topological
algebra}, Mem. Amer. Math. Soc., \textbf{7}(1951).

\bibitem{Kan} L. V. Kantorovi\v{c}, \textit{On the translocation of masses}, C. R. (Doklady) Acad. Sci. URSS (N.S.), \textbf{37}(1942), 199-201.

\bibitem{KanRub} L. V. Kantorovi\v{c}, G. \v{S}. Rubin\v{s}te\u{\i}n, \textit{On a functional space and certain
extremum problems}, Dokl. Akad. Nauk SSSR (N.S.),
\textbf{115}(1957), 1058-1061.

\bibitem{Ke1} D. Kerr, \textit{Matricial quantum Gromov-Hausdorff
distance}, J. Funct. Anal., \textbf{205}(2003), no. 1, 132-167.
arXiv:math.OA/0207282.

\bibitem{Ke2} D. Kerr, \textit{Dimension and dynamical entropy for metrized
$C^\ast$-algebras}, Comm. Math. Phys., \textbf{232}(2003), no. 3,
501-534. arXiv:math.OA/0211043.

\bibitem{KeLi} D. Kerr, H. Li, \textit{On Gromov-Hausdorff convergence for operator metric spaces},
J. Operator Theory, to appear, arXiv:math.OA/0411157.

\bibitem{Lan} G. Landi, \textit{An introduction to noncommutative spaces and their geometries}, Lecture Notes in Physics. New Series m: Monographs, 51.
Springer-Verlag, Berlin, 1997. arXiv:hep-th/9701078.

\bibitem{Li1} H. Li, \textit{Order-unit quantum Gromov-Hausdorff
distance}, J. Funct. Anal., \textbf{231}(2006), no. 2, 312-360.
arXiv:math.OA/0312001.

\bibitem{Li2} H. Li, \textit{$C^\ast$-algebraic quantum Gromov-Hausdorff distance},
arXiv:math.OA/0312003.

\bibitem{Li3} H. Li, \textit{Compact quantum metric spaces and ergodic actions of compact quantum groups}, J. Funct. Anal., to appear, arXiv:math.OA/0411178.

\bibitem{Li4} H. Li, \textit{Metric aspects of noncommutative homogeneous spaces},
arXiv:0810.4694.

\bibitem{OzRi} N. Ozawa, M. A. Rieffel, \textit{Hyperbolic group $C^\ast$-algebras and
free-product $C^\ast$-algebras as compact quantum metric spaces},
Canad. J. Math., \textbf{57}(2005), no. 5, 1056-1079.
arXiv:math.OA/0302310.

\bibitem{Pa} V. I. Paulsen, \textit{Completely bounded maps and dilations}, Pitman Research Notes in Mathematics Series, 146.
Longman Scientific \& Technical, Harlow; John Wiley \& Sons, Inc.,
New York, 1986.

\bibitem{Ri1} M. A. Rieffel, \textit{Metrics on states from actions of compact groups}, Doc.
Math., \textbf{3}(1998), 215-229. arXiv:math.OA/9807084.

\bibitem{Ri2} M. A. Rieffel, \textit{Metrics on state spaces}, Doc.
Math., \textbf{4}(1999), 559-600. arXiv:math.OA/9906151.

\bibitem{Ri5} M. A. Rieffel, \textit{Group $C^*$-algebras as compact quantum metric
spaces}, Doc. Math., \textbf{7}(2002), 605-651.
arXiv:math.OA/0205195.

\bibitem{Ri3} M. A. Rieffel, \textit{Gromov-Hausdorff distance for quantum metric
spaces}, Mem. Amer. Math. Soc. \textbf{168}(2004), no. 796, 1-65.
arXiv:math.OA/0011063.

\bibitem{Ri4} M. A. Rieffel, \textit{Matrix algebras converge to the sphere for quantum Gromov-Hausdorff distance},
Mem. Amer. Math. Soc. \textbf{168}(2004), no. 796, 67-91.
arXiv:math.OA/0108005.

\bibitem{Ri6} M. A. Rieffel, \textit{compact quantum metric spaces}, Operator algebras, quantization, and noncommutative geometry, 315-330,
Contemp. Math., 365, Amer. Math. Soc., Providence, RI, 2004.
arXiv:math.OA/0308207.

\bibitem{Ri7} M. A. Rieffel, \textit{Vector bundles and Gromov-Hausdorff distance},
arXiv:math.MG/0608266.

\bibitem{Ri8} M. A. Rieffel, \textit{Leibniz seminorms for ``Matrix algebras converge to the sphere"},
arXiv:0707.3229.

\bibitem{Val} P. Valtancoli, \textit{Projectors for the fuzzy
sphere}, Modern Phys. Lett. A, \textbf{16}(2001), no. 10, 639-645.
arXiv:hep-th/0101189.

\bibitem{Va} J. C. V\'{a}rilly, \textit{An introduction to noncommutative
geometry}, EMS Series of Lectures in Mathematics. European
Mathematical Society (EMS), Z\"{u}rich, 2006.
arXiv:physics/9709045.

\bibitem{Vo1} D. V. Voiculescu, \textit{Symmetries of some reduced free product
$C^\ast$-algebras}, Operator algebras and their connections with
topology and ergodic theory (Bu\c steni, 1983), 556-588, Lecture
Notes in Math., 1132, Springer, Berlin, 1985.

\bibitem{Vo3} D. V. Voiculescu, \textit{On the existence of quasicentral approximate units relative to normed
ideals. I}, J. Funct. Anal., \textbf{91}(1990), no. 1, 1-36.

\bibitem{VoDyNi} D. V. Voiculescu, K. J. Dykema, A. Nica, \textit{Free random
variables}, A noncommutative probability approach to free products
with applications to random matrices, operator algebras and
harmonic analysis on free groups. CRM Monograph Series, 1.
American Mathematical Society, Providence, RI, 1992.

\bibitem{Wu4} W. Wu, \textit{Locally pre-$C^\ast$-equivalent algebras}, Proc. Amer.
Math. Soc., \textbf{131}(2003), no. 2, 555-562.

\bibitem{Wu1} W. Wu, \textit{Non-commutative metric topology on
matrix state spaces}, Proc. Amer. Math. Soc. \textbf{134}(2006),
no. 2, 443-453. arXiv:math.OA/0410587.

\bibitem{Wu3} W. Wu, \textit{Quantized Gromov-Hausdorff distance}, J. Funct. Anal. \textbf{238}(2006), no. 1,
58-98. arXiv:math.OA/0503344.


\end{thebibliography}

\end{document}